\documentclass{elsart}
\usepackage{amsmath,amsfonts,amscd,amssymb,graphicx}
\journal{Physica D}
\newcommand{\trace}{\text{\rm trace }}
\newcommand{\RR}{{\mathbb R}}
\newcommand{\CC}{{\mathbb C}}
\newcommand{\Range}{\mathop\mathrm{Range}\nolimits}
\def\const{\text{\rm constant}}

\begin{document}
\begin{frontmatter}
\title{\bf An efficient shooting algorithm for Evans function calculations in large systems\thanksref{tit}}
\thanks[tit]{This research was initiated as a result of discussions during the Summer 2005 AIM workshop: Stability Criteria for Multi-Dimensional Waves and Patterns.  We thank the AIM for its hospitality and for providing a conducive environment for research, and the other participants for stimulating discussion sessions on numerical Evans function calculations in large systems. Thanks in particular to Chris Jones for the suggestion to seek an efficient representation for exterior products, which led directly to this work, and to David Bindel for a number of helpful discussions on continuation of subspaces and numerical eigenvalue problems. Numerical computations were carried out using STABLAB, an interactive MATLAB-based stability package developed by J.H.}
\author{Jeffrey Humpherys}
\address{Brigham Young University, Provo, UT 84602, USA}
\ead{jeffh@math.byu.edu}
\author{Kevin Zumbrun\thanksref{KZ}}
\address{Indiana University, Bloomington, IN 47405, USA}
\ead{kzumbrun@indiana.edu}
\thanks[KZ]{Partially supported under NSF grants number DMS-0070765 and DMS-0300487.}
\begin{abstract}
In Evans function computations of the spectra of asymptotically constant-coefficient linear operators, a basic issue is the efficient and numerically stable computation of subspaces evolving according to the associated eigenvalue ODE.  For small systems, a fast, shooting algorithm may be obtained by representing subspaces as single exterior products \cite{AS,Br.1,Br.2,BrZ,BDG}.  For large systems, however, the dimension of the exterior-product space quickly becomes prohibitive, growing as $\binom{n}{k}$, where $n$ is the dimension of the system written as a first-order ODE and $k$ (typically $\sim n/2$) is the dimension of the subspace.  We resolve this difficulty by the introduction of a simple polar coordinate algorithm representing ``pure'' (monomial) products as scalar multiples of orthonormal bases, for which the angular equation is a numerically optimized version of the continuous orthogonalization method of Drury--Davey \cite{Da,Dr} and the radial equation is evaluable by quadrature.  Notably, the polar-coordinate method preserves the important property of analyticity with respect to parameters.
\end{abstract}
\end{frontmatter}
\section{Introduction}\label{intro}

A useful tool in the study of stability of traveling waves is the Evans function, an analytic function whose zeroes correspond to the eigenvalues of the linearized operator about the wave. More generally, let $L$ be a linear
differential operator with asymptotically constant coefficients along some 
preferred spatial direction $x$, and suppose that the eigenvalue equation
\begin{equation}\label{eval}
(L-\lambda)w=0
\end{equation}
may be expressed as a first-order ODE in an appropriate phase space:
\begin{equation}\label{firstorder}
\begin{aligned}
W_x&=A(x,\lambda)W, \qquad
\lim_{x\to \pm \infty } A(x,\lambda)=A_\pm (\lambda),
\end{aligned}
\end{equation}
with $A$ analytic in $\lambda$ as a function from $\CC$ to $C^1(\RR,\CC^{n\times n})$ and the dimension $k$ of the stable subspace $S_+$ of $A_+$ and dimension $n-k$ of the unstable subspace $U_-$ of $A_-$ summing to the dimension $n$ of the entire phase space.

Then, one may construct analytic bases of solutions of \eqref{eval}, say $w_1^+, \dots, w_k^+$ and $w_{k+1}^-, \dots, w_n^-$ respectively, spanning the manifolds of solutions decaying as $x\to +\infty$ and $-\infty$ by essentially initializing them at infinity with values from the stable (resp. unstable) subspace of $A_+$ (resp. $A_-$) and solving toward $x=0$ using the ODE \eqref{firstorder}.  The Evans function is then defined as
\begin{equation}\label{evansdef}
D(\lambda):=\det\begin{pmatrix} W_1^+ & \cdots & W_k^+&W_{k+1}^- & \cdots & W_n^-\end{pmatrix}_{|x=0},
\end{equation}
where each $W_i$ is the solution of \eqref{firstorder} corresponding to $w_i$; for details, see, e.g., \cite{AGJ,KS,GZ,S.1,Z} and references therein.  Analogous to the characteristic polynomial for a finite-dimensional operator, $D(\cdot)$ is analytic in $\lambda$ with zeroes corresponding in both location and multiplicity to the eigenvalues of the linear operator $L$.

Numerical approximation of the Evans function breaks into two problems: the computation of analytic bases for stable (resp. unstable) subspaces of $A_+$ (resp. $A_-$) and the solution of ODE \eqref{firstorder} on a sufficiently large interval $x\in [0,M]$ (resp. $x\in [-M, 0]$).  In both steps, it is desirable to preserve the fundamental property of analyticity in $\lambda$, which is extremely useful in computing roots by winding number or other topological considerations.

A difficulty in the second problem is numerical stiffness for $k$, $n-k\ne 1$, due to the need to resolve modes of different exponential decay (resp. growth) rates.  This may be overcome in elegant fashion by working in the exterior product space $W_1^+\wedge \cdots \wedge W_k^+ \in \CC^{\binom{n}{k}}$ (resp. $W_{k+1}^-\wedge \cdots \wedge W_n^- \in \CC^{\binom{n}{n-k}}$), for which the desired subspace appears as a single, maximally stable (resp. unstable) mode, the Evans determinant then being recovered through the isomorphism
\[
\det\begin{pmatrix} W_1^+ & \cdots & W_k^+&W_{k+1}^- & \cdots & 
W_n^-\end{pmatrix}
\sim
W_1^+\wedge \cdots \wedge W_k^+ \wedge
W_{k+1}^-\wedge \cdots \wedge W_n^-.\\
\]
The first instance of this ``exterior-product method'' in the Evans function context seems to be a computation carried out by Pego in the Appendix of \cite{AS}. The method was subsequently rediscovered and further developed by Brin et al. \cite{Br.1,Br.2,BrZ} and, independently, by Bridges et al.\ \cite{BDG,AlB}.  See also, the earlier ``compound-matrix method'' introduced by Gilbert and Backus \cite{GB} and also Ng and Reid \cite{NR.1,NR.2,NR.3,NR.4} for the numerical solution of stiff ODE, of which it may be regarded as a coordinate-independent implementation.

The computation of an initializing analytic basis at plus (resp. minus) spatial infinity is likewise straightforward in the exterior-product framework, since it reduces to the calculation of a single eigenvector.  Two quite satisfactory approaches to this problem were given in \cite{BrZ} and \cite{BDG}, each of order $\binom{n}{k}^3$ equal to the cost of a matrix inversion or the multiplication of two matrices in dimension $\binom{n}{k}\times \binom{n}{k}$: negligible compared to the cost of integrating the exterior-product version of \eqref{firstorder}.

Together, these two steps give an extremely fast and well-conditioned shooting algorithm for the computation of the Evans function, for small values of $n$. However, for equations of large dimension $n$, such as those that arise in complicated physical systems or through transverse discretization of a multi-dimensional problem on a cylindrical domain \cite{LPSS}, the exterior-product method quickly becomes impractical, since the typical value $k\sim n/2$ leads to a working dimension $\binom{n}{n/2}$ growing as $n^{n/2}$. For example, for the typical values $n\sim 10^2$ found in \cite{LPSS}, this is clearly out of computational range.  The development of new numerical methods suitable for stability analysis of large systems was cited in a recent A.I.M. workshop on {\it Stability Criteria for Multi-Dimensional Waves and Patterns}, May 2005, as one of three overarching problems facing the traveling wave community in the next generation \cite{J}.

Discussion of this problem has so far centered mostly on boundary-value methods.  For example, one may always abandon the Evans function formulation and go back to direct discretization/Galerkin techniques, hoping to optimize perhaps by multi-pole type expansions on a problem-specific basis.  However, this ignores the useful structure, and associated dimensionality reduction, encoded by existence of the Evans function.

Alternatively, Sandstede \cite{S.2} has suggested to work within the Evans function formulation, but, in place of the high-dimensional shooting methods described above, to recast \eqref{firstorder} as a boundary-value problem with appropriate projective boundary conditions, which may be solved in the original space $\CC^n$ for individual modes by robust and highly-accurate collocation/continuation techniques.  This would reduce the cost to polynomial order $C(n)kn^2$, where $C(\cdot)$ counts the number of mesh points times evaluations per mesh required for system size $n$. He points out, further, that if one is interested only in zeroes of $D(\cdot)$ and not component subspaces,  then the cost may be reduced by a further factor of $k$ by a root-finding algorithm computing only a single candidate eigenfunction in each of the unstable (resp. stable) subspaces (the ``bordered matrix'' method as  described in \cite{BDF,Ke}).  In \cite{HSZ}, there were presented correspondingly efficient $O(n^3)$ initialization routines prescribing analytic basis/projection for the stable (resp. unstable) subspace of $A_+$ (resp. $A_-$), also in the original space $\CC^n$ without reference to exterior products.

However, up to now, it is not known how to recover analyticity of the Evans function in numerically well-conditioned fashion by the above-described collocation methods. The reason is precisely the (normally advantageous) property that errors are uniformly spatially distributed for such methods, whence the relative error near spatial infinity, where solutions exponentially decay, is prohibitively large to track dependence on initializing conditions. Likewise, the appealing simplicity/ease of coding of shooting methods is lost in this approach.

Motivated by these circumstances, we introduce in this note an alternative, shooting method designed for large systems, couched like collocation methods in the original (relatively) lower-dimensional space $\CC^n$ but preserving the useful properties of analyticity, simplicity, and good numerical conditioning enjoyed by the exterior-product method.  Indeed, being based on standard matrix operations, our method is substantially easier than the exterior-product method to code, and has an inherent parallel structure that may be exploited in a transparent fashion by the use of a numerical package incorporating parallel matrix routines; likewise, because it is carried out in minimal coordinates, there is no need to take advantage of sparse matrices such as occur in the exterior-product method.

The basis of the new method is to represent the exterior products of the columns of $W_\pm$ in ``polar coordinates'' $(\Omega, \gamma)_\pm$, where the columns of $\Omega_+ \in \CC^{n\times k}$ and $\Omega_- \in \CC^{(n-k)\times k}$ are orthonormal bases of the subspaces spanned by the columns of 
$W_+:= \begin{pmatrix} W_1^+ & \cdots & W_k^+ \end{pmatrix}$ and 
$W_-:= \begin{pmatrix} W_{k+1}^- & \cdots & W_n^- \end{pmatrix}$, $W_k^\pm$ 
defined as in \eqref{evansdef}, i.e., $W_+=\Omega_+ \alpha_+$, 
$W_-=\Omega_- \alpha_-$, and 
$ \gamma_\pm:= \det \alpha_\pm$,
so that
\[
W_1^+\wedge \cdots \wedge W_k^+ \wedge= \gamma_+ 
(\Omega_+^1\wedge \cdots \wedge \Omega_+^k),
\]
where $\Omega_+^j$ denotes the $j$th column of $\Omega_+$,
and likewise
\[
W_{k+1}^-\wedge \cdots \wedge W_n^-=\gamma_- 
(\Omega_-^1 \wedge \cdots \wedge \Omega_-^{n-k}).
\]  
The idea is that the projectivized, angular flow in $\Omega$ should remain numerically well-conditioned since orthonormality is enforced, whereas the 
radial equation, being scalar and linear, 
is automatically well-conditioned and evaluable by simple quadrature. 
Indeed, this turns out to be the case, as described in the remainder of the 
paper.
 
Just as the exterior-product method has ties to the much earlier compound-matrix method, our polar-coordinate approach has ties to another well-known method, that of ``continuous orthogonalization'', introduced by Drury \cite{Dr} and Davey \cite{Da} as an alternative to the compound-matrix technique.  Specifically, our angular equation in principle computes the same $\Omega$ computed in the continuous orthogonalization method. The new feature of our method is the radial equation, which restores the important property of analyticity with respect to parameters with no loss of numerical conditioning.  This represents a subtle but important departure in point of view, in that we truly compute exterior products and not subspaces or bases thereof, as in past interpretations of the (standard) continuous orthogonalization method \cite{B}. As discussed in Section \ref{stabilized}, this gives us considerable flexibility in choosing a numerically optimal implementation.  For the examples considered, we found excellent results (and optimal computation cost among methods considered) with the original equation of Drury \cite{Dr}; see Section \ref{numerical}.  However, this flexibility might be useful in more numerically challenging situations.

We point out that Bridges \cite{B} has developed a clever ``biorthogonal'' 
variant of Davey's continuous orthogonalization method that also preserves 
analyticity (see Remark \ref{nodamping});
indeed, not only the minors of $\Omega$ but also individual 
columns vary analytically with respect to a parameter.  However, this method 
does not preserve orthogonality but only biorthogonality 
$\tilde \Omega \Omega=I_{k\times k}$ with a simultaneously computed  
``left basis'' $\tilde \Omega$, hence, in addition to requiring twice the 
computation time due to the doubled variable $(\Omega, \tilde \Omega)$, 
appears to be inherently less stable than the other methods considered here. 
As far as we know, this method has never been implemented numerically; 
it was presented in \cite{B} as an interesting ``dual'' version of the exterior-product method.

Our numerical experiments indicate that the polar-coordinate method is quite competitive in mesh-size requirement with the exterior-product method.  Thus, the break-even dimension for the polar coordinate vs. exterior product method, taking into account competing effects of nonlinear function calls vs. higher dimension, seems to be about $n\ge 6$ ($n\ge 8$ for optimized exterior product with  sparse matrix solver, which at the moment does not exist). Detailed comparisons are given in Section \ref{numerical}.  As compared to collocation/continuation methods, we expect as for any shooting method that there is a transition size $n_*$ above which the stability advantages of collocation outweigh the speed and ease of computation of our algorithm.  In particular, for ``medium-sized'' systems such as occur in large but genuinely one-dimensional systems such as magnetohydrodynamics ($n=15$) or combustion with many species ($n\sim 10$ or $10^2$), we expect (though, so far, no study has been made for large systems with either method) that our algorithm is at least competitive with the present alternatives for root-finding.  And, for the moment, it is the only choice for calculating winding numbers in medium-sized systems or above.

Finally, we point out that, for ultra-large systems for which shooting may not be well-advised, our algorithm may equally well serve as the basis of analyticity-preserving collocation methods, since uniform errors in orthonormal bases give good tracking of subspaces along the whole real line. Thus, it seems perhaps useful in this context as well.

\section{The algorithm}\label{algorithm}
 
\subsection{Derivation}\label{derivation}
Our starting point is a comment by Chris Jones \cite{J} that the representation in the exterior-product method of products $W_1\wedge \cdots \wedge W_k$ as the direct sum of products of standard basis elements is extremely inefficient, since the less than $k \times n$-dimensional manifold of ``pure'' (i.e., monomial) products is quite small in the $\binom{n}{k}$-dimensional space of direct sums, and one should therefore seek more efficient coordinates for the computation. A natural choice is to represent exterior products $\Lambda$ in ``polar coordinates'' as $(\gamma, \Omega)$, where radius $\gamma \in \CC$ is a complex scalar and angle $\Omega \in \CC^{n\times k}$ is matrix of orthonormal column vectors $\Omega^*\Omega=I_{k\times k}$ whose columns span the subspace spanned by the factors of $\Gamma$, with $\gamma$ chosen so that the product of $\gamma$ with the exterior product of the columns of $\Omega$ is equal to $\Lambda$. This representation is unique modulo transformations $(\gamma, \Omega)\to(\gamma/\det U, \Omega U)$, $U\in \CC^{k\times k}$ is unitary.  The set of angles $\Omega$ may be recognized as the Stiefel manifold described in \cite{BDF}, a standard coordinatization of the Grasmannian manifold of linear $k$-dimensional subspaces over $\CC^n$.

In these coordinates, our computations have a concrete, linear-algebraic interpretation in which no reference to exterior products appears. Hereafter, let $\Re M:=(1/2)(M+M^*)$ denote the Hermitian part of a matrix $M$ and $\Im M:=(1/2)(M-M^*)$ the skew-Hermitian part, and $'$ denote $d/dx$. Denoting by $W^+ \in \CC^{n\times k}$ and $W^- \in \CC^{n\times n-k}$ the matrices $(W_1^+, \dots, W_k^+)$ and $(W_{k+1}^-, \dots, W_{n}^-)$ from whose columns the Evans function is determined by \eqref{evansdef}, we have
\begin{equation}\label{rel}
\begin{aligned}
W^+&= \Omega^\pm \alpha^+; \qquad
\det \alpha^+=\gamma^+,
\end{aligned}
\end{equation}
and likewise for $W^-$. Thus, \eqref{evansdef} becomes simply
\begin{equation}\label{def2}
\begin{aligned}
D(\lambda)= \gamma^+\gamma^-\det(\Omega^+, \Omega^-)_{x=0}.
\end{aligned}
\end{equation}

Fix the $\Omega$-evolution by the choice 
\begin{equation}\label{Omega}
\Omega^*\Omega'=0,
\end{equation}
removing the ambiguity in representation. This may be achieved by a unitary transformation $\Omega \to \Omega U$, $U\in \CC^{k\times k}$ satisfying $U'=-\Omega^*\Omega' U$. (Note that $-\Omega^*\Omega$ is skew-symmetric, since $0=(\Omega^* \Omega)_x=2\Re (\Omega^* \Omega_x)$.)  This choice is optimal in the sense of minimizing arc length in the Stiefel coordinates, as discussed in \cite{BDF}.  It may also be recognized as the one prescribed by a standard continuation algorithm of Kato \cite{Kat} with the orthogonal projection $P(\lambda):= \Omega \Omega^*$; see \cite{HSZ} or Remark \ref{orth} below.

Comparing $W'=AW=A\Omega \alpha$ and $W'= (\Omega \alpha)' = \Omega'\alpha + \Omega \alpha'$, we obtain
\begin{equation}\label{calc1}
\begin{aligned}
A\Omega \alpha= \Omega'\alpha + \Omega \alpha'.
\end{aligned}
\end{equation}
Multiplying on the left by $\Omega^*$ and invoking \eqref{Omega} and $\Omega^*\Omega=I$, the key equation
\begin{equation}\label{alphaeq}
\begin{aligned}
\alpha'&= (\Omega^* A \Omega) \alpha
\end{aligned}
\end{equation}
relates the two coordinatizations of the desired subspace.

By Abel's equation, we obtain from \eqref{alphaeq} and definition $\gamma:=\det \alpha$ 
\begin{equation}\label{gammaeq}
\begin{aligned}
\gamma'&= \trace (\Omega^*A\Omega)\gamma.
\end{aligned}
\end{equation}
Finally, substituting \eqref{alphaeq} into \eqref{calc1}, multiplying on the right by $\alpha^{-1}$, and rearranging, we obtain
\begin{equation}\label{Omegaeq}
\begin{aligned}
\Omega'&= (I-\Omega\Omega^*) A \Omega,
\end{aligned}
\end{equation}
completing our description of the polar-coordinate flow.
As described in the introduction, ODE \eqref{Omegaeq} is exactly
the continuous orthogonalization method of Drury \cite{Dr}.

\begin{rem}\label{abel}
Equation \eqref{gammaeq} gives another sense in which $\Omega$ is a minimal, or ``neutral'' choice of basis; namely, it is the unique choice for which the generalized Abel formula \eqref{gammaeq} holds up to complex phase. (It holds in modulus for any orthonormal basis choice.)
\end{rem}

{\bf Rescaled radial flow.}
Since we ultimately evaluate $\gamma$ at $x=0$, we may 
strategically introduce, similarly as in \cite{Br.1,Br.2,BrZ}
for the exterior-product method, the rescaled variables
\begin{equation}\label{rescaleq}
\tilde \gamma_\pm(x):=\gamma_\pm(x) e^{- \trace (\Omega^*A\Omega)_\pm x}
\end{equation}
for which the flow near $x=\pm \infty$ is asymptotically trivial.
Our complete algorithm thus becomes
\begin{equation}\label{polaralg}
\begin{aligned}
\Omega'&=(I- \Omega \Omega^*)A(x,\lambda)\Omega\\
\tilde \gamma'&= \trace (\Omega^* (A-A_-) \Omega)\tilde \gamma,
\end{aligned}
\end{equation}
with 
\begin{equation}\label{def3}
\begin{aligned}
D(\lambda)= \tilde \gamma^+\tilde \gamma^-\det(\Omega^+, \Omega^-)_{x=0}.
\end{aligned}
\end{equation}

Summarizing, we have 

\begin{prop}\label{exist}
For any choice of $\tilde \gamma(-\infty)$ and $\Omega(-\infty))$ 
with columns spanning the unstable subspace $U^+$ of 
$A_-:=\lim_{x\to -\infty} A(x,\lambda)$, there are 
unique solutions of \eqref{alphaeq}, \eqref{polaralg} such that
\begin{equation}
\begin{aligned}
\lim_{x\to -\infty}\Omega(x)&=\Omega(-\infty),\\ 
\lim_{x\to -\infty}\tilde \gamma(x)
&=\tilde \gamma(-\infty),\\
\lim_{x\to -\infty}\alpha(x)
e^{-\trace(\Omega^*A\Omega(-\infty))x}
&=\tilde \gamma(-\infty) I,\\
\end{aligned}
\end{equation}
and these satisfy $W=\Omega \alpha$, $\tilde \gamma = e^{-\trace(\Omega^*A\Omega(-\infty))x} \det \alpha$, where $W$ is a solution of $W'=AW$. 
\end{prop}

\begin{pf}
Standard asymptotic ODE theory \cite{CL} and the above calculations relating $W$, $(\gamma, \Omega)$, and $\tilde \gamma$.
\end{pf}

\subsection{Initialization at infinity}
To complete the description of our method, it remains to prescribe the initializing values $\Omega(\pm \infty)$, $\tilde \gamma (\pm \infty)$ in Proposition \ref{exist}, taking care to preserve analyticity with respect to $\lambda$. Recall the following standard result of Kato.

\begin{prop}[{\cite[\S~II.4.2]{Kat}}]\label{kato}
Let $P(\lambda)$ be an analytically varying projection on a simply connected domain $\Omega$.  Then, the linear analytic ODE
\begin{equation}\label{basisode}
r_j'=[P',P]r_j; \quad r_j(\lambda_0)=r_j^0
\end{equation}
defines a global analytically varying basis $\{r_j(\lambda)\}$ of the associated invariant subspace $\Range P(\lambda)$, where ``$\mbox{~}^\prime$'' denotes $d/d\lambda$ and $[A,B]:=AB-BA$ the commutator of $A$ and $B$. More generally,
\begin{equation}\label{kode}
S'= [P',P]S; \quad S(\lambda_0)=I.
\end{equation}
defines a globally analytic coordinate change such that
\begin{equation}\label{const1}
S^{-1}PS\equiv \const=P_0.
\end{equation}
\end{prop}

\begin{pf}
Relation \eqref{const1} follows from $(S^{-1}PS)'=0$, and may be established directly by using the key relations $PP'P=0$ and $(I-P)(I-P)'(I-P)=0$, which in turn follow by differentiation of the projective equation $P^2-P=0$. Observing that $S^{-1}$ satisfies the ``transpose'' ODE,
\begin{equation}\label{adode}
(S^{-1})'= -S^{-1}S'S^{-1}= -S^{-1}[P',P]SS^{-1}=-S^{-1}[P',P],
\end{equation}
we have that both $S$ and its inverse satisfy linear analytic ODEs, hence have global bounded analytic solutions in $\Omega$ by standard analytic ODE theory \cite{CL}.
Finally, $\Range P= S\Range P_0$ is spanned by the columns of $S R_0$, where the columns of $R_0$ are chosen to span $\Range P_0$, verifying \eqref{basisode}.
\end{pf}

\begin{rem}\label{orth}
Let $R(\lambda)=\begin{pmatrix} r_1 & \cdots & r_k\end{pmatrix}\in \CC^{n\times k}$, be the matrix of basis vectors of $\Range P$ defined by ODE \eqref{basisode} and $L\in \CC^{k\times n}$ be the matrix whose rows form the dual basis of $\Range P^*$, $LR\equiv I_k\in  \CC^{k\times k}$. Then (see Proposition 2.5, \cite{HSZ}), the flow \eqref{basisode} is uniquely determined by the property \begin{equation}\label{LRnormal}
LR'\equiv L'R\equiv 0.
\end{equation}
\end{rem}

\begin{rem}\label{unit}
From \eqref{adode} we see that $S$ (hence $R$) is unitary if $P$ is self-adjoint (i.e., orthogonal projection), since in that case $[P',P]=-[P',P]$, so that $S^{-1}$ and $S^*$ satisfy the same ODE with same initial conditions $I$.  Likewise, the relation $P=SP_0S^{-1}=SP_0S^*$ shows that $S$ is unitary only if $P$ is self-adjoint with respect to some fixed coordinate system (any coordinates for which $P_0$ is self-adjoint).
\end{rem}

Proposition \ref{kato} describes a method to generate a globally analytic matrix $W^-(\lambda)$ (the matrix $R(\lambda)$ above) with columns spanning the unstable subspace $U^-$ of $A_-(\lambda)$ and similarly for $S^+$, $A_+$. Efficient numerical implementations are described in \cite{HSZ}.

Likewise, we may efficiently compute a matrix $\Omega(-\infty, \lambda)$ at each $\lambda$ whose columns form an orthonormal basis for $U$, for example by either the SVD or the QR-decomposition.  This need not even be continuous with respect to $\lambda$.  Equating
\[
\Omega(-\infty, \lambda)\alpha(-\infty, \lambda)= W^-(\lambda),
\]
we obtain
\[
\alpha(-\infty, \lambda)= \Omega^*(-\infty, \lambda) W^-(\lambda)
\]
and therefore the exterior product of the columns of $W^-$ is equal to the exterior product of the columns of $\Omega(-\infty)$ times 
\begin{equation}\label{gamma-}
\tilde\gamma(-\infty):= 
\det \alpha(-\infty, \lambda)= \det(\Omega^*(-\infty, \lambda) W^-(\lambda)).
\end{equation}
That is, the exterior product represented by polar coordinates $(\gamma, \Omega)(-\infty, \lambda)$, with $\gamma(-\infty, \lambda)$ defined as in \eqref{gamma-}, is the same as the exterior product of the columns of $W^-(\lambda)$ appearing in the definition of the Evans function, in particular, analytic with respect to $\lambda$.

With these initializing values $(\gamma, \Omega)(-\infty, \lambda)$, we may then efficiently solve \eqref{polaralg} for $\Omega$ from $x=\pm \infty$ to  $x=0$ to obtain $\Omega_\pm(0)$ using any reasonable adaptive Runga--Kutta solver, with good numerical conditioning.  (The reason, similarly as for the exterior-product method, is that $\Omega(\pm \infty)$ is now an attractor for the flow in the direction we are integrating; see discussion, 
\cite{AGJ,GZ,Br.1,Br.2,BrZ,BDG}.)
Combining, we obtain $D(\lambda)$ through formula \eqref{def3}.

\section{Further elaborations}

\subsection{Numerical stabilization}\label{stabilized}

We briefly describe various alternatives to the basic scheme
\eqref{polaralg}, designed to improve numerical stability in
sensitive situations.
In the examples we considered, such additional stabilization
was not necessary.

{\bf Angle equation.}
It was reported soon after its 
introduction that the basic continuous orthogonalization method
\eqref{Omegaeq} can in some situations
suffer from numerical instability, by Davey \cite{Da},
who suggested as a variant the 
{\it Davey} (or generalized inverse) {\it method}
\begin{equation}
\label{davey}
\Omega'= (I-\Omega\Omega^+) A \Omega,
\end{equation}
where $\Omega^+:= (\Omega^* \Omega)^{-1}\Omega^*$ denotes the 
generalized inverse.  
The point is that $\Re(\Omega^*\Omega')=0$ for \eqref{Omegaeq}
only on the Stiefel manifold, and so level sets of the
error $D(\Omega):=\Omega^*\Omega-I$ are in general not preserved,
the error equation being
\begin{equation}\label{druryerror}
D'= - 2\, \Re ( D \Omega^* A\Omega).
\end{equation}

Davey's method, on the other hand, is derived precisely from the global
requirement $\Omega^*\Omega'=0$, and so preserves {\it all} level
sets of $D$, with associated error equation
\begin{equation}\label{daveyerror}
D'= 0.
\end{equation}
That is, \eqref{davey} corrects for spurious growth modes
of \eqref{Omegaeq} in directions normal to the Stiefel manifold.

\begin{rem}
With the corresponding modification $\tilde \gamma'= \trace (\Omega^+(A-A_\pm)\Omega)\tilde \gamma$ of \eqref{Omegaeq}, \eqref{davey} gives an alternative implementation of the full polar-coordinate method, valid on or off the Stiefel manifold.
\end{rem}

An alternative stabilization of \eqref{Omegaeq} is the {\it damped Drury method}
\begin{equation}\label{CR}
\begin{aligned}
\Omega'&=(I- \Omega \Omega^*)A(x,\lambda)\Omega
+c\Omega(I-\Omega^*\Omega), \\
\end{aligned}
\end{equation}
suggested by \cite{CR}, with $c>0$ chosen sufficiently large with respect to the matrix norm $\|A\|$ of $A$.  This evidently agrees with \eqref{gammaeq}--\eqref{Omegaeq} on the manifold $\mathcal{G}:=\{\Omega: \, \Omega^*\Omega\equiv I\}$, but, thanks to the new penalty term on the right-hand side of the $\Omega$ equation, has the additional favorable property that $\mathcal{G}$ is not only invariant but attracting under the flow.  For, defining $D(\Omega)=\Omega^*\Omega-I$ as above, we have the error equation
\begin{equation}\label{dampeddruryerror}
D'=-2 c(I+D)D - 2\, \Re ( D \Omega^* A\Omega).
\end{equation}
Observing that $\Omega$ and $D$ have matrix norm one, we obtain for $c> \|A\|$ that
\[
D'\le -2 cD - 2 \, \Re( D \Omega^* A\Omega)< 0.
\]

\begin{rem}\label{frobenius}
A brief calculation reveals that the stabilizing term $\Omega(I-\Omega^*\Omega)$ is the steepest descent/gradient flow for the Frobenius norm squared of $D$. In this sense, it is the optimal, least-squares correction, corresponding approximately with orthogonal projection onto the Stiefel manifold $D=0$.  
\end{rem}

Though \eqref{CR} in principle remains stable in numerically sensitive regimes where \eqref{davey} may fail, in practice, the large value of $c$ that is required for stability introduces numerical stiffness imposing unreasonable restrictions on mesh size; see Section
\ref{numerical}.
A more attractive alternative in this situation is the {\it damped Davey method}
\begin{equation}
\label{dampeddavey}
\Omega'= (I-\Omega\Omega^+) A \Omega
+c\Omega(I-\Omega^*\Omega) \\
\end{equation}
obtained by combining the above two stabilization techniques, which is globally exponentially attracting for any $c>0$, with error equation
\begin{equation}\label{dampeddaveyerror}
D'= -2 c(I+D)D. 
\end{equation}
{This would be our own suggestion in numerically sensitive regimes.}

Alternatively, one might employ the {\it damped Bridges--Reich method}
\begin{equation}
\label{BrRe}
\Omega'=-2\Im[(I- \Omega \Omega^*)A(x,\lambda)]\Omega+c \Omega(I-\Omega^*\Omega), 
\end{equation}
proposed in \cite{BrRe}, for which the Stiefel manifold $\mathcal{G}$ likewise is attracting for any $c>0$. This agrees with \eqref{Omegaeq} on $\mathcal{G}$,  since
\[
[(I- \Omega \Omega^*)A]^*\Omega = 
A^*(I- \Omega \Omega^*)\Omega=0;
\]
moreover, the {\it undamped Bridges--Reich method} obtained by setting $c=0$ has a skew-Hermitian form favorable for geometric integrators preserving Lie group structure \cite{BrRe}.  On the other hand, off of the Stiefel manifold, these methods do not preserve the desired subspace, since they are not of the canonical form 
\begin{equation}\label{canonical}
\Omega'=A\Omega + \Omega B,
\end{equation}
where $B=\alpha'\alpha^{-1}$, obtained by multiplying \eqref{calc1} on the right by $\alpha^{-1}$; thus, they are in some sense trading errors in the normal direction for new errors tangential to the Stiefel manifold.

Finally, we mention the {\it polar factorization method}, a discrete orthogonalization method introduced by Higham \cite{Hig}, based on the observation that the factor $\Omega$ in the polar factorization $W= \Omega H$ of a matrix $W$, $\Omega$ orthonormal and $H$ symmetric, is the closest point (in usual Frobenius matrix norm) to $W$ on the Stiefel manifold.  In this split-step method, single steps of an explicit scheme approximating the original, linear flow $W'=AW$ are alternated with projection via polar-factorization onto the Stiefel manifold.  Though we did not test this scheme, Higham \cite{Hig} reported superior performance in numerically sensitive situations, as compared to geometric integrator-based continuous orthogonalization methods.  For further discussion of the continuous orthogonalization method and its variants, see \cite{AlB,B,BrRe} and references therein.

\begin{rem}\label{higrem}
In applying Higham's method to Evans function computations, we expect that it is important to substitute for $W'=AW$ the asymptotically trivial rescaled ODE, $W'=(A-\trace(\Omega^* A \Omega(\pm\infty))I)W$, similarly as for the exterior-product method \cite{Br.1,Br.2,BrZ} or the radial equation in \eqref{polaralg}. In this approach, $\tilde \gamma$ is obtained as the product $\Pi_j \det H_j$ of determinants of symmetric factors $H$ in the projection steps $W=\Omega H \to \Omega$.  Alternatively, one could integrate Drury's ODE \eqref{Omegaeq} for the ODE step and track $\tilde \gamma$ by the usual, continuous ODE 
\eqref{polaralg}.
\end{rem}

\begin{rem}\label{nodamping} An analytic variant of Davey's method is the
{bi-orthogonal method}
\begin{equation}\label{Banalytic}
\begin{aligned}
 \Omega'&=(I-  \Omega(\tilde \Omega^* \Omega)^{-1}\tilde
\Omega^*)A \Omega,\\
 \tilde \Omega'&=-(I- \tilde \Omega(\Omega^* \tilde \Omega)^{-1}
\Omega^*)A^* \tilde \Omega,\\
\end{aligned}
\end{equation}
introduced by Bridges \cite{B}.  Here, the matrix $\Omega$, and 
likewise $\tilde \Omega^*$, is analytic along with its $k\times k$ minors. 
This scheme has the property of Davey's method that level sets 
of $\tilde \Omega^* \Omega$ are preserved; in particular, 
$\tilde \Omega^*\Omega\equiv I_{k\times k}$ is an invariant manifold.  
However, analyticity of $\Omega$ and $\tilde \Omega^*$ is incompatible 
with $\tilde \Omega=\Omega$, and so $\tilde \Omega^* \Omega=I_{k\times k}$ 
does not enforce good conditioning of $\Omega$.
Apparently, the requirement of analyticity of individual columns is overly 
rigid for purposes of numerical stabilization; the advantage of the polar 
coordinate method is the flexibility to choose optimally the angular 
evolution equation.
\end{rem}

In our numerical experiments, the best and fastest performance was exhibited by the original, undamped Drury method \eqref{Omegaeq}; see Section \ref{numerical}. Indeed, our results indicate that the examples considered lie in the numerically insensitive regime for which normal instabilities remain small. However, the stability issues discussed above may be important for other, more numerically taxing problems.

{\bf Radial equation.}
Likewise,
since they are scalar quantities, we could more stably solve for radial variables $\tilde\gamma_\pm(0)$ by quadrature, using formulae
\begin{equation}\label{gammaeval-}
\begin{aligned}
\tilde \gamma_-(0)&=e^{\int_{-\infty}^0 \trace(\Omega^*A\Omega(x)-
\Omega_-^*A_-\Omega_-) dx}\tilde \gamma(-\infty),\\
\tilde\gamma_+(0)&=e^{\int_0^{+\infty} -\trace(\Omega^*A\Omega(x)-
\Omega_+^*A_+\Omega_+) dx}\tilde\gamma(+\infty).\\
\end{aligned}
\end{equation}
However, in practice, this does not seem necessary, as the $\Omega$-equation appears to be the limiting factor determining accuracy/allowable mesh size in \eqref{polaralg}.
{\it The rescaling \eqref{rescaleq} on the other hand is critical for good numerical performance, as is the analogous rescaling for the exterior-product method \cite{Br.1,Br.2,BrZ}, giving a speedup on the order of the spectral radius of $A(x,\lambda)$.}

\subsection{Continuous initialization}

Though the initializing step at plus and minus spatial infinity is a one-time cost, hence essentially negligible, we point out that this step too may be carried out more efficiently by an evolution scheme in $\lambda$ parallel to the one carried out in $x$ in Section \ref{derivation}.  Defining $W_-:=W^-(-\infty, \lambda)$, $\gamma_-:=\gamma(-\infty, \lambda)$, and $Q:=[P',P]$, recall that the $\lambda$ evolution of $W_-$ is again given by a linear ODE,
\[
W_-'=Q(\lambda)W_-. 
\]

Thus, we may apply exactly the same steps as in the previous section to obtain a well-conditioned $C^k$ evolution scheme for $(\Omega_-, \tilde \gamma_-)$ of
\begin{equation}\label{lambdacont}
\begin{aligned}
\Omega_-'&=(I- \Omega_- \Omega_-^*)Q(\lambda)\Omega_-
+c(\|Q\|)\Omega_-(I-\Omega_-^*\Omega_-), \\
\tilde \gamma_-'&= \trace (\Omega_-^* Q \Omega_-)\tilde \gamma_-,
\end{aligned}
\end{equation}
initializing $(\tilde \gamma_{-0}, \Omega_{-0})= (1, \Omega_{-0})$ at some base point $\lambda_0$.

With this modification, equations \eqref{gammaeval-} simplify to
\begin{equation}\label{sgammaeval}
\begin{aligned}
\tilde \gamma_-(0)&=e^{\int_{-\infty}^0 \trace(\Omega^*A\Omega(x)-
\Omega_-^*A_-\Omega_-) dx}\tilde \gamma(-\infty),\\
\tilde \gamma_+(0)&=e^{\int_0^{+\infty} -\trace(\Omega^*A\Omega(x)-
\Omega_+^*A_+\Omega_+) dx}\tilde \gamma(+\infty).\\
\end{aligned}
\end{equation}

Note that the above calculation gives the interesting information of the winding number of $\tilde \gamma_-$ over one circuit around a closed $\lambda$-contour, which is not necessarily zero, or even an integer.

\begin{rem}\label{colloctrack}
A similar, but more complicated scheme could be used to restore analyticity in general collocation methods, by further tracking in $x$ the values of $\tilde \gamma$ relating the bases obtained by Kato's algorithm applied in variable $\lambda$ with respect to orthogonal projection;
for, the associated $x$-evolution depends only on the associated subspace, which can be well-approximated, rather than exterior product or directions of individual solutions, which cannot.
Combined with the above computation relating analytic bases at $x=\infty$ to those obtained through orthogonal projection, we obtain an analytic scheme for which the only information required is knowledge (to reasonable tolerance) of subspaces at each $x$.  Of course, as pointed out in the introduction, a simpler solution would be to use a collocation scheme based on the algorithm of this paper, for which no such corrections would be necessary.
\end{rem}

\section{Numerical comparisons}\label{numerical}

In this section, we compare our new method for Evans function computation with the exterior-product method described in Section \ref{intro}, and afterward, briefly, with various alternative continuous orthogonalization methods substituted in the angular equation.

As a test system, we consider solitary waves of the ``good'' Boussinesq equation
\begin{equation}
\label{boussinesq_equation}
u_{t t} = u_{x x} - u_{x x x x} - (u^2)_{x x},
\end{equation}
which have the form ($\xi=x-s t$)
\begin{equation}
\label{boussinesq_profile}
\bar{u}(\xi) = \frac{3}{2} (1-s^2) \mbox{sech}^2\left(\frac{\sqrt{1-s^2}}{2}\xi\right),
\end{equation}
where the wave speed $s$ satisfies $|s|<1$.  We remark that this is the same system studied in \cite{AS}, and is known to be stable when $\frac{1}{2}\leq |s|<1$ and unstable when $|s|<1/2$.

By linearizing \eqref{boussinesq_equation} about the traveling wave \eqref{boussinesq_profile}, we arrive at the eigenvalue problem
\begin{equation}
\lambda^2 u - 2 s \lambda u' = (1-s^2) u'' - u'''' - (2 \bar{u}u)'',
\end{equation}
which can be written as a first-order system \eqref{firstorder}, where
\begin{equation}
\label{a}
A(x,\lambda) \begin{pmatrix}0&1&0&0\\0&0&1&0\\0&0&0&1\\-\lambda^2-2\bar{u}_{x x}&2\lambda s-4\bar{u}_x&(1-s^2)-2\bar{u}&0\end{pmatrix}.
\end{equation}
For $\lambda$ in the right-hand plane, we have that \eqref{a} spectrally separates into two growth and two decay modes, that is, $k=2$.  We remark that this model is a great test-case for Evans function computation since it captures one of the chief numerical obstacles, that is, overcoming multi-mode growth and decay in the left and right subspaces, respectively.  It is precisely this difficulty that motivated the development and use of the exterior-product and continuous orthogonalization methods in the first place.

\subsection{Algorithms}

Using the exterior-product method, we lift $A(x,\lambda)$ into
exterior-product space to get
\begin{equation}\label{a2}
A^{(2)}(x,\lambda) \begin{pmatrix}
0&1&0&0&0&0\\
0&0&1&1&0&0\\
2\lambda s-4\bar{u}_x&(1-s^2)-2\bar{u}&0&0&1&0\\
0&0&0&0&1&0\\
\lambda^2+2\bar{u}_{x x}&0&0&(1-s^2)-2\bar{u}&0&1\\
0&\lambda^2+2\bar{u}_{x x}&0&-2\lambda
s+4\bar{u}_x&0&\phantom{111}0\phantom{111}
\end{pmatrix}.
\end{equation}
To maintain analyticity, we use Kato's method (Proposition \ref{kato}, see also \cite{BrZ,BDG}) for analytically choosing the (simple) eigenvectors $r_+(-M)$ and $r_-(M)$, which correspond to the largest growth and decay modes, respectively, at the numerical approximates of negative and positive infinity, $\pm M$ (we used $M=8$), where the eigenprojection $P$ is obtained as $(l^* r)^{-1} r l^{*}$ for any left and right eigenvectors $l$ and $r$ for the eigenvalue of \eqref{a2} of smallest (resp. largest) real part, obtained through a standard matrix routine. To integrate \eqref{basisode} we use Euler's first-order method for convenience.  We then evolve these vectors from $x=\pm M$ to $x=0$ using a standard numerical ODE solver (RKF45) and compute the Evans function via the wedge product at $x=0$ as in \eqref{evansdef}.

With our new method, we similarly evolve, analytically in $\lambda$, the eigenvectors at the numerical end states $\pm M$ so that we can determine the $\mbox{det}(\Omega^* W^\pm)$ multipliers in \eqref{gammaeval-}. We likewise use Kato's method, where the eigenprojection $P$ is obtained as $(L^*R)^{-1}RL^*$ for any matrices $L$ and $R$ with columns forming orthonormal bases for the left and right stable (resp. unstable) subspace of $A$, obtained by the singular-value decomposition (SVD).  For a more efficient, but slightly more complicated algorithm, see \cite{HSZ}.

We likewise determine orthonormal bases $\Omega(\pm \infty)$ at each $\lambda$ for the stable (resp. unstable) subspace of $A_+$ (resp. $A_-$) via the SVD and initialize $\tilde \gamma(\pm M)$ through \eqref{gamma-}.  Finally, we evolve $(\Omega, \tilde \gamma)$ from $x=\pm M$ to $x=0$ (via RKF45) and compute the Evans function as the Grammian determinant \eqref{def2}, as described in Section \ref{algorithm}.

\begin{figure}[t]
\centerline{
\mbox{\includegraphics[width=3.00in]{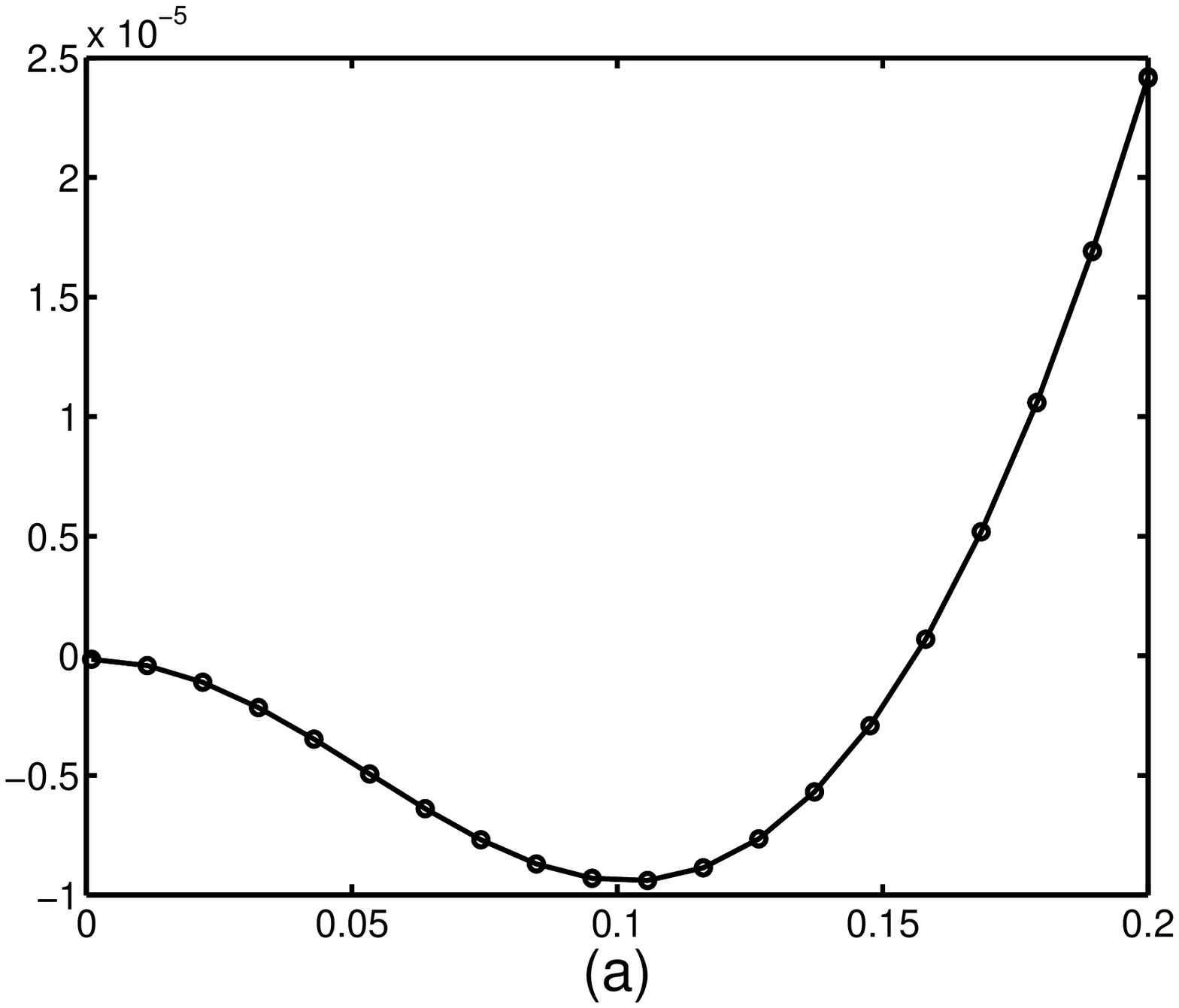}}
\mbox{\includegraphics[width=2.90in]{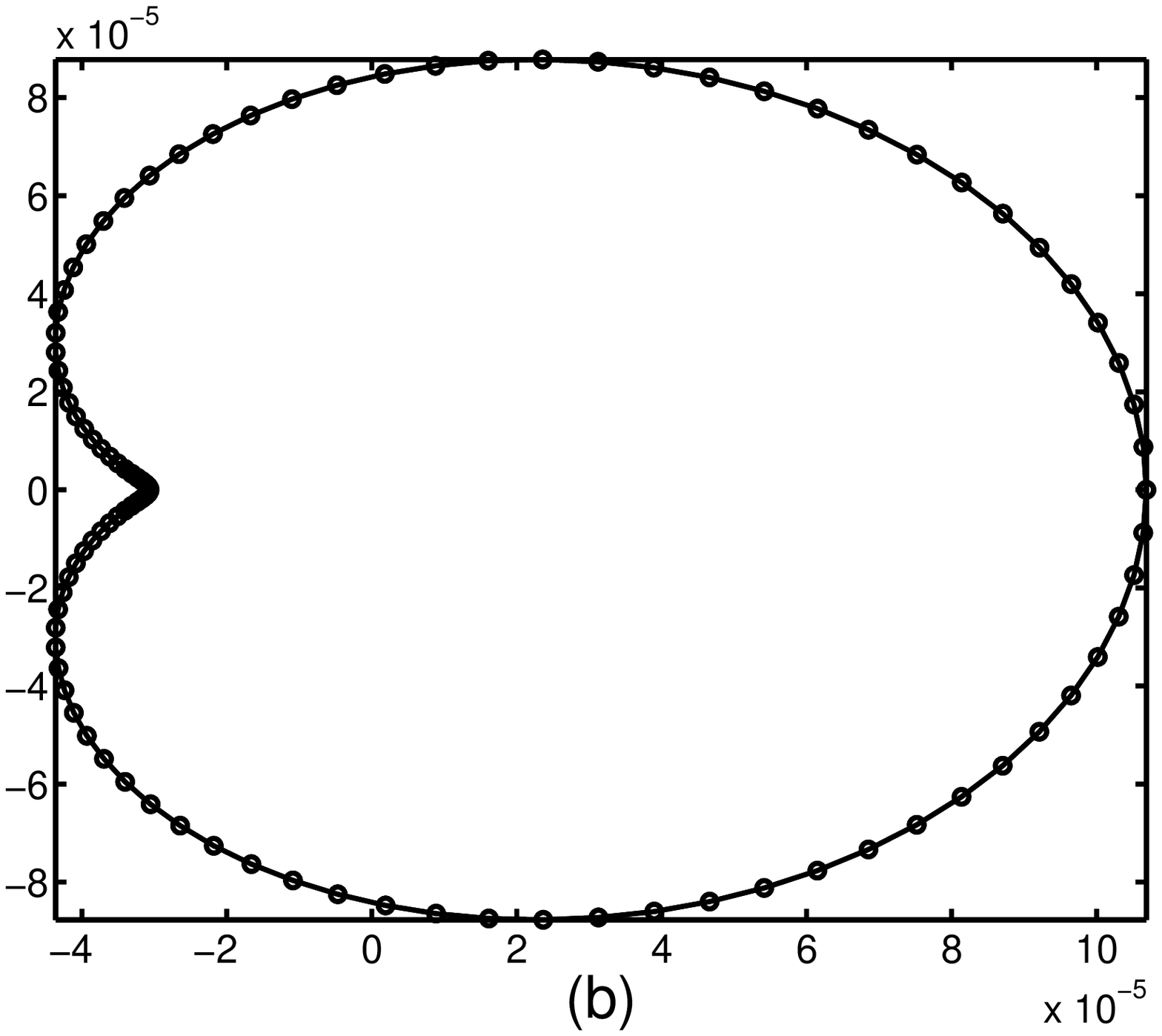}}
}
\caption{We evaluate the Evans function of the ``good'' Boussinesq system for the unstable pulse having wave speed $s=0.4$.  We compute $(a)$ the graph of the Evans function along real axis from $\lambda=0$ to $\lambda=0.2$ and $(b)$ the image of the closed contour $\Gamma(t) = 0.16 + 0.05 e^{2\pi i t}$, where $0\leq t < 1$.}
\label{boussinesq_fig}
\end{figure}

\subsection{Results}

In Figure \ref{boussinesq_fig}, we compute the Evans function for the unstable pulse with $s=0.4$.  We perform two Evans function computations
for each contour, using our new method and the exterior-product method. We then compare results.   In Figure \ref{boussinesq_fig}(a), we compute the graph of the Evans function along the real axis from $\lambda=0$ to $\lambda=0.2$.  From this we see that the graph crosses through zero near $\lambda=0.155$, indicating an unstable eigenvalue there.  We remark that the graphs for both methods were plotted, however since they are virtually
indistinguishable, it only appears as though one curve is present. In Figure \ref{boussinesq_fig}(b), we compute the Evans function about the contour $\Gamma(t) = 0.16 + 0.05 e^{2\pi i t}$, where $0\leq t < 1$, and plot its image.  The interior of this contour contains the above mentioned unstable eigenvalue, and so we expect the origin to be contained in the interior of the image, as it is.  This is a second verification of the unstable eigenvalue.  In this second test, both methods are likewise graphed and overlap to the point that they are also indistinguishable.  Indeed the two graphs differ by an absolute difference of $1.4\times 10^{-9}$ and a relative difference of $4.6\times 10^{-5}$.  In both old and new shooting algorithms, the absolute and relative tolerances are set to $10^{-8}$ and $10^{-6}$, respectively.

\subsection{Performance}

Function evaluation for \eqref{polaralg} requires as many as five matrix-matrix multiplications, whereas the exterior-product method requires a single matrix-vector multiplication.  Because of this, our new method is actually slower for the ``good'' Boussinesq, which is a relatively small system.  To leading order, the operational count for our method grows as $2 k n^2 + 3 k^2 n$.  By contrast, the exterior-product method grows as $\binom{n}{k}^2$, and is faster than our method when $n=4$ and $k=2$.  Indeed for $k\sim n/2$, we expect the break-even point to be at around $n=6$.

We remark that $A^{(k)}(x,\lambda)$ becomes sparse as $n$ gets large ($k\sim n/2$).  We can show that the number of non-zero entries (and hence operations for a sparse matrix-vector evaluation) is exactly $(k(n-k)+1)\binom{n}{k}$.  Even though a sparse-matrix package would offer substantial computational savings over a non-sparse one, the size is still prohibitive for large systems.  In fact, with a sparse improvement, the break-even point for function evaluation relative to our method is only extended only to around $n=8$.

\begin{table}
\begin{center}
\begin{tabular}{| l c | c | c | c | c | c |}\hline
Method & $c$ & $\|D\|^2$ & Mesh & Time & Abs. & Rel. \\\hline
Damped Drury \eqref{CR} & 0 & 1.6(-10) & 54 & 5.4 & 2.0(-11) & 6.5(-09)\\
& 1 & 4.1(-11) & 54 & 5.8 & 2.0(-11) & 6.3(-09)\\
& 5 & 2.8(-11) & 56 & 5.9 & 2.4(-11) & 7.7(-09)\\
& 10 & 2.0(-07) & 77 & 8.3 & 7.3(-10) & 2.4(-07)\\
& 20 & 2.8(-07) & 124 & 13.4 & 1.4(-09) & 4.6(-07)\\
& 1600 & 3.9(-07) & 7727 & 846 & 1.8(-09) & 5.7(-07)\\\hline
Damped Davey \eqref{dampeddavey}& 0 & 4.8(-10) & 54 & 6.7 & 9.2(-10) & 3.0(-07)\\
& 1 & 1.0(-10) & 54 & 7.1 & 1.9(-11) & 6.2(-09)\\
& 5 & 1.0(-11) & 54 & 7.1 & 1.9(-11) & 6.1(-09)\\
& 10 & 5.8(-11) & 54 & 7.1 & 1.9(-11) & 6.2(-09)\\
& 20 & 3.4(-07) & 94 & 12.3 & 8.8(-10) & 2.9(-07)\\
& 1600 & 1.5(-07) & 7696 & 1051 & 1.1(-09) & 3.4(-07)\\\hline
Bridges-Reich \eqref{BrRe} & 0 & 4.8(-10) & 54 & 5.8 & 7.9(-09) & 2.5(-06)\\
& 1 & 1.1(-10) & 54 & 6.1 & 1.8(-09) & 5.8(-07)\\
& 5 & 8.6(-12) & 55 & 6.1 & 3.6(-10) & 1.2(-07)\\
& 10 & 2.4(-11) & 57 & 6.5 & 4.4(-11) & 1.4(-08)\\
& 20 & 1.2(-09) & 97 & 11.0 & 2.3(-10) & 7.4(-08)\\
& 1600 & 3.3(-08) & 7705 & 893 & 8.4(-09) & 2.7(-06)\\\hline
\end{tabular}
\caption{A comparision of methods: The third column ``Mesh'' correpsonds to the (typical) number of mesh points integrating from $x=\pm M$ to $x=0$.  The ``Time'' column measures the run-time for computing the 20 Evans function values along the contour $\gamma(t) = 0.16 + 40 i + .15 e^{2\pi i t}$.  The last two columns measure the absolute and relative difference compared to the values returned using exterior-product method to high accuracy.}
\label{comparison}
\end{center}
\end{table}

\subsection{Other Methods}

We compare the different continuous orthogonalization methods discussed in Section \ref{stabilized} and examine overall performance and accuracy.  Specifically, we measure $(i)$ how closely the trajectories ``stick'' to the Stiefel manifold by computing the distance (in Frobenius norm) of $\Omega$ at $x=0$ to the Stiefel manifold, $(ii)$ the number of mesh points that are needed for our adaptive ODE solver to maintain tolerance ($\text{AbsTol=1e-8}$ and $\text{RelTol=1e-6}$), $(iii)$ the run-time needed to compute the Evans function, and finally $(iv)$ the absolute and relative errors compared with the exterior-product method taken to high accuracy.  See Table \ref{comparison} for the results.

According to the data, the overall best method for our problem is the original undamped Drury method \eqref{polaralg}.  It is only slightly less accurate than the damped Davey method \eqref{dampeddavey} with ($c=5$) but is about 25\% faster as fewer operations are needed.  We remark that all three methods become both less accurate and (considerably) more time consuming, due to numerical stiffness, as the damping constant $c$ gets large.

We also remark that, while the Bridges-Reich method \eqref{BrRe} did not perform as well as either the Drury or Davey methods, the undamped method was actually designed for use in geometric integrators, and the damping term was only introduced to facilitate less sophisticated numerical integration methods.  Indeed the use of geometric integrators seems like a very good direction to explore for Evans function computation, and we intend on exploring this issue in 
future work.


\section{Concluding remarks}\label{conclusion}

The numerical results show that the above-described polar-coordinate algorithm is indeed feasible, and compares favorably in performance even for low-dimensional systems to the exterior-product method that is the current standard. Due to better dimensional scaling, it should work well also for large systems, whereas the exterior-product method quickly becomes dimensionally infeasible. As test problems of medium size, we intend next to investigate stability of traveling waves in magnetohydrodynamics and detonation with large number of reactant species ($n\sim 15$). A longer term project might be to implement polar-coordinate based collocation methods for extremely large systems.

It seems worth emphasizing a philosophical point associated with the new algorithm that is simple but possibly of wider use, concerning the apparently conflicting goals of numerical well-conditioning vs. maintaining analyticity.  Namely, in the context of exterior products, an optimal strategy is to projectivize, choosing the most numerically convenient basis for the associated subspace (Grasmannian), then {\it correct} for the resulting loss of analyticity by an appropriate scalar factor.  Because scalar equations are always numerically well-conditioned, the final step costs essentially nothing. That is, we may recover analytic continuation of subspaces by a simple, well-conditioned post-processing step appended to a more standard $C^k$ continuation routine.

Finally, we point out that our experiments indicate that essentially any existing continuous orthogonalization methods should suffice for stability computations in low- and medium-frequency regimes relevant for sectorial operators- in particular, those arising in the reaction--diffusion context of \cite{LPSS}.  For high-frequency computations in numerically sensitive situations, it may be necessary to use the more stable schemes described in Section \ref{stabilized}. Here, we may draw on the resource of the very active community studying numerical continuous orthogonalization.
\bibliographystyle{abbrv}
\bibliography{evans}

\end{document}